\documentclass{svmult}
\usepackage[all]{xy}
\usepackage{amssymb}
\usepackage{amsopn}
\usepackage{amsmath}

\newcommand{\cf}{{\em cf.}\ }

\newcommand{\ko}{\: , \;}

\newcommand{\opname}[1]{\operatorname{\mathsf{#1}}}

\newcommand{\rep}{\opname{rep}\nolimits}

\renewcommand{\rep}{\opname{rep}\nolimits}

\newcommand{\Z}{\mathbb{Z}}
\newcommand{\N}{\mathbb{N}}
\newcommand{\Q}{\mathbb{Q}}
\newcommand{\C}{\mathbb{C}}

\newcommand{\iso}{\stackrel{_\sim}{\rightarrow}}

\newcommand{\Hom}{\opname{Hom}}

\newcommand{\Ext}{\opname{Ext}}

\newcommand{\ca}{{\mathcal A}}

\newcommand{\cc}{{\mathcal C}}
\newcommand{\cd}{{\mathcal D}}

\begin{document}

\title*{Categorification of acyclic cluster algebras: an introduction}
\author{Bernhard Keller}
\institute{Universit\'e Paris 7 -- Denis Diderot\\
    UFR de Math\'ematiques\\
   Institut de Math\'ematiques de Jussieu, UMR 7586 du CNRS \\
   Case 7012\\
   2, place Jussieu\\
   75251 Paris Cedex 05\\
   France }

\maketitle

\centerline{\emph{To Murray Gerstenhaber and Jim Stasheff}}
\vskip1cm

\begin{abstract}
This is a concise introduction to Fomin-Zelevinsky's cluster algebras
and their links with the representation theory of quivers in
the acyclic case. We review the definition cluster algebras (geometric,
without coefficients),
construct the cluster category and present the bijection between
cluster variables and rigid indecomposable objects of the cluster
category.
\end{abstract}


\date{version du 21/10/2007, modifi\'ee le 21/10/2007}

\section{Introduction}

\subsection{Context}
Cluster algebras were invented by S.~Fomin and A.~Zelevinsky
\cite{FominZelevinsky02} in the spring of the year 2000 in a project
whose aim it was to develop a combinatorial approach to the results obtained
by G.~Lusztig concerning total positivity in algebraic groups
\cite{Lusztig96}
on the one hand and canonical bases in quantum groups
\cite{Lusztig90} on the
other hand (let us stress that canonical bases were discovered
independently and simultaneously by M.~Kashiwara \cite{Kashiwara90}).
Despite great
progress during the last few years
\cite{FominZelevinsky03} \cite{BerensteinFominZelevinsky05}
\cite{FominZelevinsky07},
we are still relatively
far from these initial aims. Presently, the best results on the link between
cluster algebras and canonical bases are probably those of
C.~Geiss, B.~Leclerc and J.~Schr\"oer
\cite{GeissLeclercSchroeer05}
\cite{GeissLeclercSchroeer06}
\cite{GeissLeclercSchroeer06a}
but even
they cannot construct canonical bases from cluster variables
for the moment. Despite these difficulties, the theory of
cluster algebras has witnessed spectacular growth thanks
notably to the many links that have been discovered with
a wide range of subjects including
\begin{itemize}
\item Poisson geometry
\cite{GekhtmanShapiroVainshtein03}
\cite{GekhtmanShapiroVainshtein05} \ldots ,
\item integrable systems \cite{FominZelevinsky03b} \ldots,
\item higher Teichm\"uller spaces
\cite{FockGoncharov03}
\cite{FockGoncharov05}
\cite{FockGoncharov07a}
\cite{FockGoncharov07b} \ldots ,
\item combinatorics and the study of combinatorial
polyhedra like the Stasheff associahedra
\cite{ChapotonFominZelevinsky02}
\cite{Chapoton04}
\cite{Krattenthaler06}
\cite{FominReading05}
\cite{Musiker07}
\cite{FominShapiroThurston06}
\ldots ,
\item commutative and non commutative algebraic geometry,
in particular the study of stability
conditions in the sense of Bridgeland \cite{Bridgeland02}
\cite{Bridgeland06}, Calabi-Yau algebras
\cite{Ginzburg06},
Donaldson-Thomas invariants \cite{Szendroi07}
\cite{Kontsevich07a} \cite{Kontsevich07} \cite{KontsevichSoibelman07}
\ldots,
\item and last not least the representation theory of quivers and
finite-dimensional algebras, cf. for example the surveys
\cite{BuanMarsh06} \cite{Ringel07} \cite{Reiten06}.
\end{itemize}
We refer to the introductory
papers \cite{Zelevinsky02} \cite{FominZelevinsky03a}
\cite{Zelevinsky04} \cite{Zelevinsky05} \cite{Zelevinsky07}
and to the cluster algebras portal \cite{Fomin07}
for more information on cluster algebras
and their links with other parts of mathematics.

The link between cluster algebras and quiver representations follows
the spirit of categorification: One tries to interpret cluster algebras
as combinatorial (perhaps $K$-theoretic) invariants associated with
categories of representations. Thanks to the rich structure of
these categories, one can then hope to prove results
on cluster algebras which seem beyond the scope of the purely combinatorial
methods. It turns out that the link becomes especially beautiful
if we use a {\em triangulated category} constructed from the category
of quiver representations, the so-called cluster category.

In this brief survey, we will review the definition of cluster
algebras and Fomin-Zelevinsky's classification theorem for
cluster-finite cluster algebras \cite{FominZelevinsky03}. We will then
recall some basic notions on the representations of a quiver without
oriented cycles, introduce the cluster category and describe its link
with the cluster algebra.

\section{Cluster algebras}

The cluster algebras we will be interested in are associated with antisymmetric
matrices with integer coefficients. Instead of using matrices, we will use
quivers (without loops and $2$-cycles), since they are easy to visualize and well-suited
to our later purposes.

\subsection{Quivers}
Let us recall that a \emph{quiver} $Q$ is an oriented graph. Thus, it
is a quadruple given by a set $Q_0$ (the set of vertices), a set $Q_1$
(the set of arrows) and two maps $s:Q_1 \to Q_0$ and $t:Q_1\to Q_0$
which take an arrow to its source respectively its target. Our quivers
are `abstract graphs' but in practice we draw them as in this example:
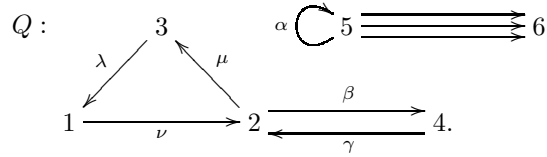
\[ Q:
\xymatrix{ & 3 \ar[ld]_\lambda & & 5 \ar@(dl,ul)[]^\alpha \ar@<1ex>[rr] \ar[rr] \ar@<-1ex>[rr] & & 6 \\
  1 \ar[rr]_\nu & & 2 \ar@<1ex>[rr]^\beta \ar[ul]_\mu & & 4.
  \ar@<1ex>[ll]^\gamma }
\]
A \emph{loop} in a quiver $Q$ is an arrow $\alpha$ whose source
coincides with its target; a \emph{$2$-cycle} is a pair of distinct
arrows $\beta\neq\gamma$ such that the source of $\beta$ equals the
target of $\gamma$ and vice versa. It is clear how to define
\emph{$3$-cycles}, \emph{connected components} \ldots . A quiver is
\emph{finite} if both, its set of vertices and its set of arrows, are
finite.

\subsection{Seeds and mutations}

Fix an integer $n\geq 1$.
A \emph{seed} is a pair $(R,u)$, where
\begin{itemize}
\item[a)] $R$ is a finite quiver without loops or $2$-cycles with vertex
set $\{1, \ldots, n\}$;
\item[b)] $u$ is a free generating set $\{u_1, \ldots, u_n\}$ of the field
$\Q(x_1, \ldots, x_n)$ of fractions of the polynomial ring $\Q[x_1, \ldots, x_n]$
in $n$ indeterminates.
\end{itemize}
Notice that in the quiver $R$ of a seed, all arrows between any two
given vertices point in the same direction (since $R$ does not have
$2$-cycles).  Let $(R,u)$ be a seed  and $k$ a vertex of $R$. The
\emph{mutation} $\mu_k(R,u)$ of $(R,u)$ at $k$ is the seed $(R',u')$,
where
\begin{itemize}
\item[a)] $R'$ is obtained from $R$ as follows:
\begin{itemize}
\item[1)] reverse all arrows incident with $k$;
\item[2)] for all vertices $i\neq j$ distinct from $k$, modify the number
of arrows between $i$ and $j$ as follows:
\[
\begin{array}{|c|c|}\hline
R & R'  \\ \hline
\xymatrix@=0.3cm{i \ar[rr]^{r} \ar[rd]_{s} & &  j \\
& k \ar[ru]_t & } &
\xymatrix@=0.3cm{i \ar[rr]^{r+st} & & j  \ar[ld]^t\\
& k \ar[lu]^{s}} \\\hline
\xymatrix@=0.3cm{i \ar[rr]^{r}  & &  j \ar[ld]^{t}\\
& k \ar[lu]^{s} & } &
\xymatrix@=0.3cm{i \ar[rr]^{r-st} \ar[rd]_{s} & &  j \\
& k \ar[ru]_t & } \\\hline
\end{array}
\]
where $r,s,t$ are non negative integers,
an arrow $\xymatrix@=0.3cm{i\ar[r]^l & j}$ with $l\geq 0$
means that $l$ arrows go from $i$ to $j$ and an arrow
$\xymatrix@=0.3cm{i\ar[r]^l & j}$ with $l\leq 0$
means that $-l$ arrows go from $j$ to $i$.
\end{itemize}
b) $u'$ is obtained from $u$ by replacing the element $u_k$ with
\begin{equation} \label{eq:exchange}
u_k'=\frac{1}{u_k} \left( \prod_{\mbox{\scriptsize arrows $i\to k$}} u_i + \prod_{\mbox{\scriptsize arrows $k\to j$}} u_j\right).
\end{equation}
\end{itemize}
In the {\em exchange relation}~(\ref{eq:exchange}),
if there are no arrows from $i$ with target $k$, the product is taken over
the empty set and equals $1$. It is not hard to see that $\mu_k(R,u)$ is
indeed a seed and that $\mu_k$ is an involution: we have $\mu_k(\mu_k(R,u))=(R,u)$.

\subsection{Examples of mutations}
Let $R$ be the
cyclic quiver
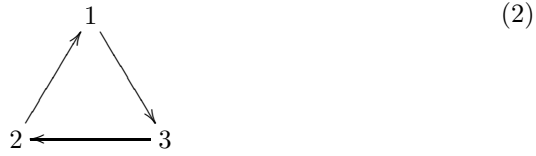
\begin{equation} \label{quiver1}
\begin{xy} 0;<0.3pt,0pt>:<0pt,-0.3pt>::
(94,0) *+{1} ="0",
(0,156) *+{2} ="1",
(188,156) *+{3} ="2",
"1", {\ar"0"},
"0", {\ar"2"},
"2", {\ar"1"},
\end{xy}
\end{equation}
and $u=\{x_1, x_2, x_3\}$. If we mutate at $k=1$, we obtain the quiver
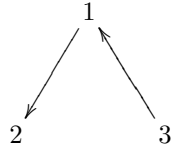
\[
\begin{xy} 0;<0.3pt,0pt>:<0pt,-0.3pt>::
(92,0) *+{1} ="0",
(0,155) *+{2} ="1",
(188,155) *+{3} ="2",
"0", {\ar"1"},
"2", {\ar"0"},
\end{xy}
\]
and the set of fractions given by $u'_1=(x_2+x_3)/x_1$, $u'_2=u_2=x_2$ and $u'_3=u_3=x_3$.
Now, if we mutate again at $1$, we obtain the original seed. This is a
general fact: Mutation at $k$ is an involution. If, on the other hand, we
mutate $(R', u')$ at $2$, we obtain the quiver
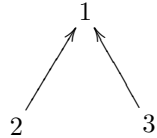
\[
\begin{xy} 0;<0.3pt,0pt>:<0pt,-0.3pt>::
(87,0) *+{1} ="0",
(0,145) *+{2} ="1",
(167,141) *+{3} ="2",
"1", {\ar"0"},
"2", {\ar"0"},
\end{xy}
\]
and the set $u''$ given by $u''_1=u'_1=(x_2+x_3)/x_1$,
$u'_2=\frac{x_1 + x_2 + x_3}{x_1 x_2}$ and
$u''_3=u'_3=x_3$.

Let us consider the following, more complicated quiver
glued together from four $3$-cycles:
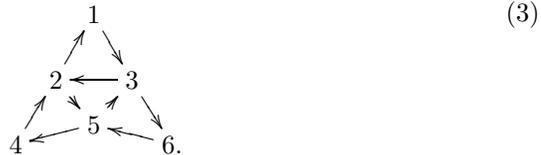
\begin{equation} \label{quiver2}
\begin{xy} 0;<0.4pt,0pt>:<0pt,-0.4pt>::
(74,0) *+{1} ="0",
(38,62) *+{2} ="1",
(110,62) *+{3} ="2",
(0,123) *+{4} ="3",
(74,104) *+{5} ="4",
(148,123) *+{6.} ="5",
"1", {\ar"0"},
"0", {\ar"2"},
"2", {\ar"1"},
"3", {\ar"1"},
"1", {\ar"4"},
"4", {\ar"2"},
"2", {\ar"5"},
"4", {\ar"3"},
"5", {\ar"4"},
\end{xy}
\end{equation}
If we successively perform mutations at the vertices $5$, $3$, $1$ and $6$,
we obtain the sequence of quivers (we use \cite{KellerQuiverMutationApplet})
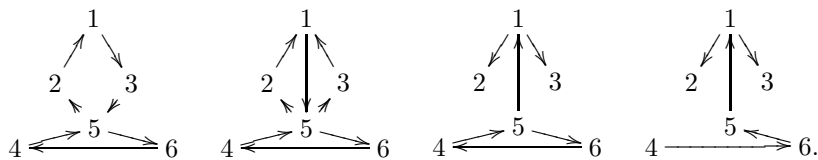
\[
\quad
\begin{xy} 0;<0.4pt,0pt>:<0pt,-0.4pt>::
(74,0) *+{1} ="0",
(38,62) *+{2} ="1",
(110,62) *+{3} ="2",
(0,123) *+{4} ="3",
(75,104) *+{5} ="4",
(148,123) *+{6} ="5",
"1", {\ar"0"},
"0", {\ar"2"},
"4", {\ar"1"},
"2", {\ar"4"},
"3", {\ar"4"},
"5", {\ar"3"},
"4", {\ar"5"},
\end{xy}
\quad
\begin{xy} 0;<0.4pt,0pt>:<0pt,-0.4pt>::
(75,0) *+{1} ="0",
(38,62) *+{2} ="1",
(110,61) *+{3} ="2",
(0,123) *+{4} ="3",
(75,104) *+{5} ="4",
(148,123) *+{6} ="5",
"1", {\ar"0"},
"2", {\ar"0"},
"0", {\ar"4"},
"4", {\ar"1"},
"4", {\ar"2"},
"3", {\ar"4"},
"5", {\ar"3"},
"4", {\ar"5"},
\end{xy}
\quad
\begin{xy} 0;<0.4pt,0pt>:<0pt,-0.4pt>::
(75,0) *+{1} ="0",
(38,61) *+{2} ="1",
(110,60) *+{3} ="2",
(0,122) *+{4} ="3",
(75,103) *+{5} ="4",
(148,122) *+{6} ="5",
"0", {\ar"1"},
"0", {\ar"2"},
"4", {\ar"0"},
"3", {\ar"4"},
"5", {\ar"3"},
"4", {\ar"5"},
\end{xy}
\quad
\begin{xy} 0;<0.4pt,0pt>:<0pt,-0.4pt>::
(75,0) *+{1} ="0",
(38,61) *+{2} ="1",
(110,60) *+{3} ="2",
(0,122) *+{4} ="3",
(75,103) *+{5} ="4",
(149,121) *+{6.} ="5",
"0", {\ar"1"},
"0", {\ar"2"},
"4", {\ar"0"},
"3", {\ar"5"},
"5", {\ar"4"},
\end{xy}
\quad
\]
Notice that the last quiver no longer has any oriented cycles
and is in fact an orientation of the Dynkin diagram of type $D_6$.
The sequence of new fractions appearing in these steps is
\begin{eqnarray*}
u'_5 & = &\frac{x_3 x_4 + x_2 x_6}{x_5} \ko \quad
u'_3  = \frac{x_3 x_4 + x_1 x_5 + x_2 x_6}{x_3 x_5}\ko \\
u'_1 & = & \frac{x_2 x_3 x_4 + x_3^2 x_4 + x_1 x_2 x_5 + x_2^2 x_6 + x_2 x_3x_6}{x_1 x_3 x_5} \ko\quad
u'_6=\frac{x_3 x_4 + x_4 x_5 + x_2 x_6}{x_5 x_6}\;.
\end{eqnarray*}
It is remarkable that all the denominators appearing here are monomials
and that all the coefficients in the numerators are positive.

Finally,
let us consider the quiver
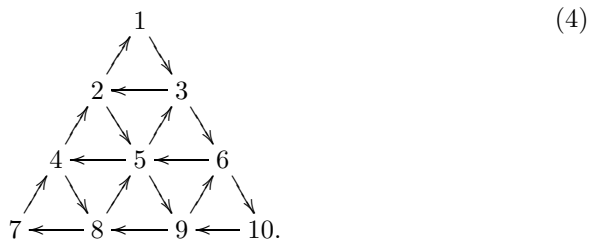
\begin{equation} \label{quiver3}
\begin{xy} 0;<0.6pt,0pt>:<0pt,-0.6pt>::
(79,0) *+{1} ="0",
(52,44) *+{2} ="1",
(105,44) *+{3} ="2",
(26,88) *+{4} ="3",
(79,88) *+{5} ="4",
(131,88) *+{6} ="5",
(0,132) *+{7} ="6",
(52,132) *+{8} ="7",
(105,132) *+{9} ="8",
(157,132) *+{10.} ="9",
"1", {\ar"0"},
"0", {\ar"2"},
"2", {\ar"1"},
"3", {\ar"1"},
"1", {\ar"4"},
"4", {\ar"2"},
"2", {\ar"5"},
"4", {\ar"3"},
"6", {\ar"3"},
"3", {\ar"7"},
"5", {\ar"4"},
"7", {\ar"4"},
"4", {\ar"8"},
"8", {\ar"5"},
"5", {\ar"9"},
"7", {\ar"6"},
"8", {\ar"7"},
"9", {\ar"8"},
\end{xy}
\end{equation}
One can show \cite{KellerReiten06} that it is impossible to transform it
into a quiver without
oriented cycles by a finite sequence of mutations. However, its mutation
class (the set of all quivers obtained from it by iterated mutations)
contains many quivers with just one oriented cycle, for example
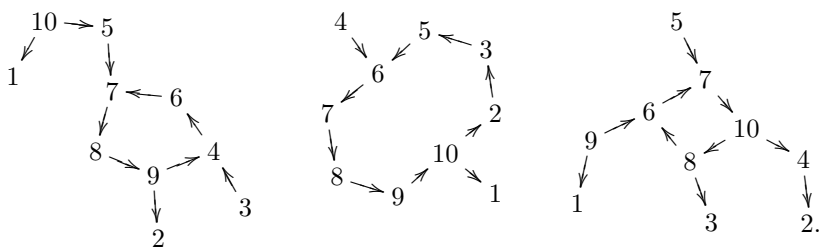
\[
\begin{xy} 0;<0.3pt,0pt>:<0pt,-0.3pt>::
(0,70) *+{1} ="0",
(183,274) *+{2} ="1",
(293,235) *+{3} ="2",
(253,164) *+{4} ="3",
(119,8) *+{5} ="4",
(206,96) *+{6} ="5",
(125,88) *+{7} ="6",
(104,164) *+{8} ="7",
(177,194) *+{9} ="8",
(39,0) *+{10} ="9",
"9", {\ar"0"},
"8", {\ar"1"},
"2", {\ar"3"},
"3", {\ar"5"},
"8", {\ar"3"},
"4", {\ar"6"},
"9", {\ar"4"},
"5", {\ar"6"},
"6", {\ar"7"},
"7", {\ar"8"},
\end{xy}
\quad\quad
\begin{xy} 0;<0.3pt,0pt>:<0pt,-0.3pt>::
(212,217) *+{1} ="0",
(212,116) *+{2} ="1",
(200,36) *+{3} ="2",
(17,0) *+{4} ="3",
(123,11) *+{5} ="4",
(64,66) *+{6} ="5",
(0,116) *+{7} ="6",
(12,196) *+{8} ="7",
(89,221) *+{9} ="8",
(149,166) *+{10} ="9",
"9", {\ar"0"},
"1", {\ar"2"},
"9", {\ar"1"},
"2", {\ar"4"},
"3", {\ar"5"},
"4", {\ar"5"},
"5", {\ar"6"},
"6", {\ar"7"},
"7", {\ar"8"},
"8", {\ar"9"},
\end{xy}
\quad\quad
\begin{xy} 0;<0.3pt,0pt>:<0pt,-0.3pt>::
(0,230) *+{1} ="0",
(294,255) *+{2.} ="1",
(169,253) *+{3} ="2",
(285,174) *+{4} ="3",
(125,0) *+{5} ="4",
(90,114) *+{6} ="5",
(161,73) *+{7} ="6",
(142,177) *+{8} ="7",
(17,150) *+{9} ="8",
(213,135) *+{10} ="9",
"8", {\ar"0"},
"3", {\ar"1"},
"7", {\ar"2"},
"9", {\ar"3"},
"4", {\ar"6"},
"5", {\ar"6"},
"7", {\ar"5"},
"8", {\ar"5"},
"6", {\ar"9"},
"9", {\ar"7"},
\end{xy}
\]
In fact, in this example, the mutation class is finite and it can be
completely computed
using, for example, \cite{KellerQuiverMutationApplet}:
It consists of $5739$
quivers up to isomorphism. The above quivers are members of the mutation
class containing relatively few arrows. The initial quiver is the unique
member of its mutation class with the largest number of arrows. Here
are some other quivers in the mutation class with a relatively
large number of arrows:
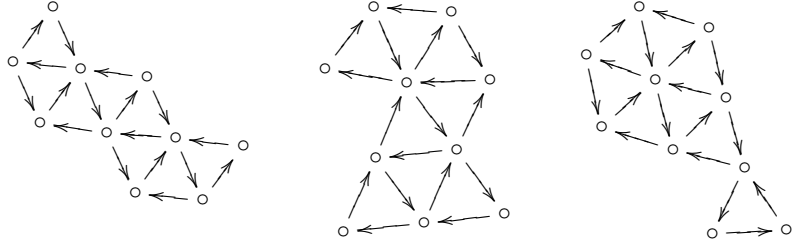
\[
\begin{xy} 0;<0.3pt,0pt>:<0pt,-0.3pt>::
(290,176) *+{\circ} ="0",
(154,235) *+{\circ} ="1",
(34,147) *+{\circ} ="2",
(50,0) *+{\circ} ="3",
(239,244) *+{\circ} ="4",
(0,69) *+{\circ} ="5",
(169,89) *+{\circ} ="6",
(205,165) *+{\circ} ="7",
(85,78) *+{\circ} ="8",
(118,159) *+{\circ} ="9",
"4", {\ar"0"},
"0", {\ar"7"},
"4", {\ar"1"},
"1", {\ar"7"},
"9", {\ar"1"},
"5", {\ar"2"},
"2", {\ar"8"},
"9", {\ar"2"},
"5", {\ar"3"},
"3", {\ar"8"},
"7", {\ar"4"},
"8", {\ar"5"},
"6", {\ar"7"},
"6", {\ar"8"},
"9", {\ar"6"},
"7", {\ar"9"},
"8", {\ar"9"},
\end{xy}
\quad
\quad
\begin{xy} 0;<0.3pt,0pt>:<0pt,-0.3pt>::
(0,78) *+{\circ} ="0",
(226,262) *+{\circ} ="1",
(23,284) *+{\circ} ="2",
(61,0) *+{\circ} ="3",
(208,92) *+{\circ} ="4",
(159,7) *+{\circ} ="5",
(125,273) *+{\circ} ="6",
(64,191) *+{\circ} ="7",
(166,180) *+{\circ} ="8",
(103,96) *+{\circ} ="9",
"0", {\ar"3"},
"9", {\ar"0"},
"1", {\ar"6"},
"8", {\ar"1"},
"6", {\ar"2"},
"2", {\ar"7"},
"5", {\ar"3"},
"3", {\ar"9"},
"5", {\ar"4"},
"8", {\ar"4"},
"4", {\ar"9"},
"9", {\ar"5"},
"7", {\ar"6"},
"6", {\ar"8"},
"8", {\ar"7"},
"7", {\ar"9"},
"9", {\ar"8"},
\end{xy}
\quad
\quad
\begin{xy} 0;<0.3pt,0pt>:<0pt,-0.3pt>::
(159,287) *+{\circ} ="0",
(252,281) *+{\circ} ="1",
(19,152) *+{\circ} ="2",
(67,0) *+{\circ} ="3",
(0,61) *+{\circ} ="4",
(200,203) *+{\circ} ="5",
(109,180) *+{\circ} ="6",
(155,26) *+{\circ} ="7",
(176,115) *+{\circ} ="8",
(87,92) *+{\circ} ="9",
"0", {\ar"1"},
"5", {\ar"0"},
"1", {\ar"5"},
"4", {\ar"2"},
"6", {\ar"2"},
"2", {\ar"9"},
"4", {\ar"3"},
"7", {\ar"3"},
"3", {\ar"9"},
"9", {\ar"4"},
"5", {\ar"6"},
"8", {\ar"5"},
"6", {\ar"8"},
"9", {\ar"6"},
"7", {\ar"8"},
"9", {\ar"7"},
"8", {\ar"9"},
\end{xy}
\]
Only $84$ among the $5739$ quivers in the mutation class contain
double arrows (and none contain arrows of multiplicity $\geq 3$).
Here is a typical example
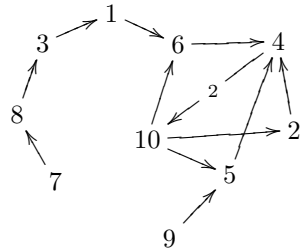
\[
\begin{xy} 0;<0.4pt,0pt>:<0pt,-0.4pt>::
(89,0) *+{1} ="0", (262,111) *+{2} ="1", (24,29) *+{3} ="2",
(247,27) *+{4} ="3", (201,153) *+{5} ="4", (152,30) *+{6} ="5",
(36,159) *+{7} ="6", (0,96) *+{8} ="7", (144,213) *+{9} ="8",
(123,120) *+{10} ="9", "2", {\ar"0"}, "0", {\ar"5"}, "1", {\ar"3"},
"9", {\ar"1"}, "7", {\ar"2"}, "4", {\ar"3"}, "5", {\ar"3"}, "3",
{\ar|*+{\scriptstyle 2}"9"}, "8", {\ar"4"}, "9", {\ar"4"}, "9",
{\ar"5"}, "6", {\ar"7"},
\end{xy}
\]
The quivers (\ref{quiver1}), (\ref{quiver2}) and (\ref{quiver3}) are part
of a family which appears in the study of the cluster
algebra structure on the coordinate algebra of the
subgroup of upper unitriangular matrices in 
$SL_n(\C)$, \cf \cite{GeissLeclercSchroeer06}.
The study of coordinate algebras on varieties associated
with reductive algebraic groups
(in particular, double Bruhat cells) has provided a major impetus
for the development of cluster algebras, \cf
\cite{BerensteinFominZelevinsky05}.

\subsection{Definition of cluster algebras}
Let $Q$ be a finite quiver without loops or $2$-cycles with
vertex set $\{1, \ldots, n\}$. Consider the seed $(Q,x)$
consisting of $Q$ and the set $x$ formed by the variables
$x_1, \ldots, x_n$. Following \cite{FominZelevinsky02} we
define
\begin{itemize}
\item the \emph{clusters with respect to $Q$} to be the
sets $u$ appearing in seeds $(R,u)$ obtained from $(Q,x)$ by
iterated mutation,
\item the \emph{cluster variables} for $Q$ to be the elements of all clusters,
\item the \emph{cluster algebra $\ca_Q$} to be the $\Q$-subalgebra of
the field $\Q(x_1, \ldots, x_n)$ generated by all the cluster variables.
\end{itemize}
Thus the cluster algebra consists of all $\Q$-linear combinations
of monomials in the cluster variables. It is useful to define
another combinatorial object associated with this recursive construction:
The \emph{exchange graph} associated with $Q$ is the graph whose
vertices are the seeds modulo simultaneous renumbering of the vertices and
the associated cluster variables and whose edges correspond to mutations.

\subsection{The example $A_3$}
Let us consider the quiver
\[
Q: \xymatrix{1 \ar[r] & 2 \ar[r] & 3}
\]
obtained by endowing the Dynkin diagram $A_3$ with a linear orientation.
By applying the recursive construction to the initial seed $(Q,x)$
one finds exactly fourteen seeds (modulo simultaneous renumbering of
vertices and cluster variables). These are the vertices of the exchange
graph, which is isomorphic to the third Stasheff associahedron \cite{Stasheff63}
\cite{ChapotonFominZelevinsky02}:
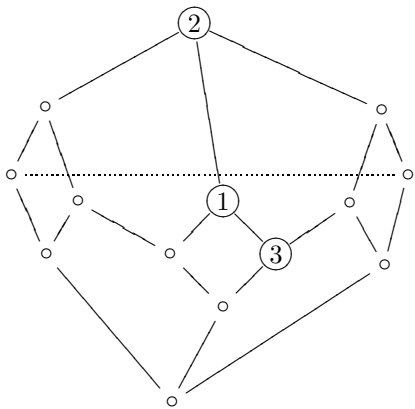
\[
\begin{xy} 0;<0.4pt,0pt>:<0pt,-0.4pt>::
(173,0) *+<8pt>[o][F]{2} ="0",
(0,143) *+{\circ} ="1",
(63,168) *+{\circ} ="2",
(150,218) *+{\circ} ="3",
(250,218) *+<8pt>[o][F]{3} ="4",
(375,143) *+{\circ} ="5",
(350,82) *+{\circ} ="6",
(152,358) *+{\circ} ="7",
(200,168) *+<8pt>[o][F]{1} ="8",
(200,268) *+{\circ} ="9",
(32,79) *+{\circ} ="10",
(33,218) *+{\circ} ="11",
(320,170) *+{\circ} ="12",
(353,228) *+{\circ} ="13",
"0", {\ar@{-}"6"},
"0", {\ar@{-}"8"},
"0", {\ar@{-}"10"},
"1", {\ar@{.}"5"},
"1", {\ar@{-}"10"},
"11", {\ar@{-}"1"},
"2", {\ar@{-}"3"},
"10", {\ar@{-}"2"},
"2", {\ar@{-}"11"},
"3", {\ar@{-}"8"},
"9", {\ar@{-}"3"},
"8", {\ar@{-}"4"},
"4", {\ar@{-}"9"},
"4", {\ar@{-}"12"},
"6", {\ar@{-}"5"},
"5", {\ar@{-}"13"},
"12", {\ar@{-}"6"},
"9", {\ar@{-}"7"},
"11", {\ar@{-}"7"},
"13", {\ar@{-}"7"},
"13", {\ar@{-}"12"},
\end{xy}
\]
The vertex labeled $1$ corresponds to $(Q,x)$, the vertex $2$ to
$\mu_2(Q,x)$, which is given by
\[
\xymatrix{1 \ar@/^1pc/[rr] & 2 \ar[l] & 3 \ar[l]} \ko \{ x_1, \frac{x_1+x_3}{x_2}, x_3\} \ko
\]
and the vertex $3$ to $\mu_1(Q,x)$, which is given by
\[
\xymatrix{1 & 2 \ar[l] \ar[r] & 3} \ko \{\frac{1+x_3}{x_1}, x_2, x_3\}.
\]
We find a total of $9$ cluster variables, namely
\begin{align*}
& x_1 \ko x_2 \ko x_3,
\frac{1+x_2}{x_1}\ko \frac{x_1+x_3+x_2x_3}{x_1x_2} \ko \frac{x_1+x_1 x_2 +x_3 + x_2 x_3}{x_1 x_2 x_3} \ko\\
& \frac{x_1+x_3}{x_2}\ko
\frac{x_1+x_1x_2+x_3}{x_2 x_3}\ko
\frac{1+x_2}{x_3}\;.
\end{align*}
Again we observe that all denominators are monomials.
Notice also that $9=3+6$ and that $3$ is the rank of the root system associated
with $A_3$ and $6$ its number of positive roots. Moreover, if we look at
the denominators of the non trivial cluster variables (those other than $x_1$, $x_2$, $x_3$), we
see a natural bijection with the positive roots
\[
\alpha_1, \alpha_1+\alpha_2, \alpha_1+\alpha_2+\alpha_3,
\alpha_2, \alpha_2+\alpha_3, \alpha_3
\]
of the root system of $A_3$, where $\alpha_1$, $\alpha_2$, $\alpha_3$ denote
the three simple roots.

\subsection{Cluster algebras with finitely many cluster variables}
The phenomena observed in the above example are explained by the following
key theorem:

\begin{theorem}[Fomin-Zelevinsky \cite{FominZelevinsky03}] Let $Q$ be
  a finite connected quiver without loops or $2$-cycles with vertex
  set $\{1, \ldots, n\}$. Let $\ca_Q$ be the associated cluster
  algebra.
\begin{itemize}
\item[a)] All cluster variables are Laurent polynomials, i.e. their
  denominators are monomials.
\item[b)] The number of cluster variables is finite iff $Q$ is
  mutation equivalent to an orientation of a simply
  laced Dynkin diagram $\Delta$. In this case, $\Delta$ is unique and
  the non trivial cluster variables are in bijection with the positive
  roots of $\Delta$; namely, if we denote the simple roots by
  $\alpha_1, \ldots, \alpha_n$, then for each positive root $\sum d_i
  \alpha_i$, there is a unique non trivial cluster variable whose
  denominator is $\prod x_i^{d_i}$.
\end{itemize}
\end{theorem}

\section{Categorification}

We refer to the books \cite{Ringel84} \cite{GabrielRoiter92}
\cite{AuslanderReitenSmaloe95} and \cite{AssemSimsonSkowronski06} for
a wealth of information on the representation theory of quivers and
finite-dimensional algebras. Here, we will only need very basic notions.

Let $Q$ be a finite quiver without oriented cycles. For example, $Q$ can
be an orientation of a simply laced Dynkin diagram or
the quiver
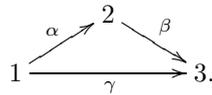
\[
\xymatrix@R=10pt{ & 2 \ar[rd]^\beta & \\
1 \ar[rr]_{\gamma} \ar[ru]^{\alpha} & & 3.
}
\]
Let $k$ be an algebraically closed field. Recall that a {\em representation
of $Q$} is a diagram of finite-dimensional vector spaces of the shape
given by $Q$. Thus, in the above example, a representation of $Q$ is
a (not necessarily commutative) diagram
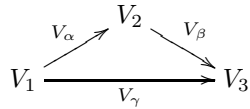
\[
\xymatrix@R=10pt{ & V_2 \ar[rd]^{V_\beta} & \\
V_1 \ar[rr]_{V_{\gamma}} \ar[ru]^{V_\alpha} & & V_3
}
\]
formed by three finite-dimensional vector spaces and three linear maps.
A \emph{morphism of representations} is a morphism of diagrams. We thus
obtain the \emph{category of representations} $\rep(Q)$. Notice that
it is an abelian category (since it is a category of diagrams in
an abelian category, that of finite-dimensional vector spaces):
Sums, kernels and cokernels in the category $\rep(Q)$ are
computed componentwise. We denote by $\cd_Q$ its bounded
derived category. Thus, the objects of $\cd_Q$ are
the bounded complexes
\[
\xymatrix{ \ldots \ar[r] & V^p \ar[r]^d & V^{p+1} \ar[r] & \ldots}
\]
of representations and its morphisms are obtained from
morphisms of complexes by formally inverting all quasi-isomorphisms
(=morphisms inducing isomorphisms in homology). The category $\cd_Q$
is still an additive category (direct sums are given by direct
sums of complexes) but it is almost never abelian. In fact, it
is abelian if and only if $Q$ does not have any arrows. But it
is always {\em triangulated}. This means that $\cd_Q$ is additive
and endowed with
\begin{itemize}
\item[a)] a {\em suspension functor} $\Sigma: \cd_Q \iso \cd_Q$, namely
the functor taking a complex $V$ to $V[1]$, where $V[1]^p=V^{p+1}$
for all $p\in\Z$ and $d_{V[1]}=-d_V$;
\item[b)] a class of {\em triangles}, namely the sequences
\[
\xymatrix{U \ar[r] & V \ar[r] & W \ar[r] & \Sigma U}
\]
which are `induced' from short exact sequences of complexes.
\end{itemize}
The triangulated category $\cd_Q$ admits a \emph{Serre functor}, i.e.
an autoequivalence $S: \cd_Q \iso \cd_Q$ which makes the Serre duality
formula true: We have
\[
D\Hom(X,Y) \iso \Hom(Y,SX)
\]
bifunctorially in $X$, $Y$ belonging to $\cd_Q$, where $D$ denotes
the duality functor $\Hom_k(?,k)$ over the ground field $k$.
The {\em cluster category} is defined as the \emph{orbit category}
\[
\cc_Q = \cd_Q / (S^{-1}\circ \Sigma^2)^\Z
\]
of $\cd_Q$ under the action of the cyclic group generated by the
automorphism $S^{-1}\circ \Sigma^2$. Thus, its objects are the
same as those of $\cd_Q$ and its morphisms are defined by
\[
\Hom_{\cc_Q}(X,Y) = \bigoplus_{p\in \Z} \Hom_{\cd_Q}(X, (S^{-1}\circ \Sigma^2)^p Y).
\]
One can show \cite{Keller05} that $\cc_Q$ admits a structure of triangulated
category such that the projection functor $\cd_Q \to \cc_Q$ becomes
a triangle functor (in general, the orbit category of a triangulated
category under the action of an automorphism group is no longer triangulated).
It is not hard to see that the cluster category has finite-dimensional
morphism spaces, and that it admits a Serre functor induced by that of
the derived category. The definition of the cluster category then immediately
yields an isomorphism
\[
S \iso \Sigma^2
\]
and this means that $\cc_Q$ is $2$-Calabi-Yau: A $k$-linear triangulated
category with finite-dimensional morphism spaces is $d$-Calabi-Yau 
if it admits a Serre functor $S$ and if $S$ is isomorphic to 
$\Sigma^d$ (the $d$th power of the suspension functor) as a triangle functor.
The definition of the cluster
category is due to Buan-Marsh-Reineke-Reiten-Todorov
\cite{BuanMarshReinekeReitenTodorov04} (for arbitrary $Q$ without
oriented cycles)
and, independently and with a very different, more geometric description,
to Caldero-Chapoton-Schiffler \cite{CalderoChapotonSchiffler06}
(for $Q$ of type $A_n$).

To state the close relationship between the cluster category $\cc_Q$ and
the cluster algebra $\ca_Q$, we need some notation: For two objects $L$ and $M$
of $\cc_Q$, we write
\[
\Ext^1(L,M) = \Hom_{\cc_Q}(L,\Sigma M).
\]
Notice that it follows from the Calabi-Yau property that we have a canonical
isomorphism
\[
\Ext^1(L,M) \iso D \Ext^1(M,L).
\]
An object $L$ of $\cc_Q$ is \emph{rigid} if we have $\Ext^1(L,L)=0$. It is
\emph{indecomposable} if it is non zero and in each decomposition
$L=L_1\oplus L_2$, we have $L_1=0$ or $L_2=0$.

\begin{theorem}[\cite{CalderoKeller06}] \label{thm:2}
Let $Q$ be a finite quiver without oriented cycles with
vertex set $\{1, \ldots, n\}$.
\begin{itemize}
\item[a)] There is an explicit bijection $L \mapsto X_L$ from the set
of isomorphism classes of rigid indecomposables of the cluster category $\cc_Q$
onto the set of cluster variables of the cluster algebra $\ca_Q$.
\item[b)] Under this bijection, the clusters correspond exactly to
the \emph{cluster-tilting subsets}, i.e. the sets $T_1, \ldots, T_n$
of rigid indecomposables such that
\[
Ext^1(T_i,T_j)=0
\]
for all $i,j$.
\item[c)] If $L$ and $M$ are rigid indecomposables such that the
space $\Ext^1(L,M)$ is one-dimensional, then we have the generalized
exchange relation
\begin{equation} \label{eq:gen-exchange}
X_L = \frac{X_B + X_{B'}}{X_M}
\end{equation}
where $B$ and $B'$ are the middle terms of `the' non split triangles
\[
\xymatrix{L \ar[r] & B \ar[r] & M \ar[r] & \Sigma L} \mbox{ and }
\xymatrix{M \ar[r] & B' \ar[r] & L \ar[r] & \Sigma M}
\]
and we define
\[
X_B= \prod_{i=1}^s X_{B_i} \ko
\]
where $B=B_1 \oplus \cdots \oplus B_s$ is a decomposition into
indecomposables.
\end{itemize}
\end{theorem}

The relation~(\ref{eq:gen-exchange}) in part c) of the theorem can be generalized
to the case where the extension group is of higher dimension,
cf. \cite{CalderoKeller05a} \cite{Hubery06} \cite{XiaoXu07}.
One can show using \cite{BuanMarshReiten04} that  relation~(\ref{eq:gen-exchange})
generalizes the exchange relation~(\ref{eq:exchange})
which appeared in the definition of the mutation.

The proof of the theorem builds on work by many authors
notably
Buan-Marsh-Reiten-Todorov \cite{BuanMarshReitenTodorov07},
Buan-Marsh-Reiten \cite{BuanMarshReiten04b},
Buan-Marsh-Reineke-Reiten-Todorov \cite{BuanMarshReinekeReitenTodorov04},
Marsh-Reineke-Zelevinsky \cite{MarshReinekeZelevinsky03},
 \ldots\  and especially on Caldero-Chapoton's explicit
formula for $X_L$ proved in \cite{CalderoChapoton06} for orientations
of simply laced Dynkin diagrams. We include the formula below. Another
crucial ingredient of the proof is the Calabi-Yau property
of the cluster category.
An alternative proof of part c) was given by A.~Hubery \cite{Hubery06}
for quivers whose underlying graph is an extended simply
laced Dynkin diagram.

The theorem does shed new light on cluster
algebras. In particular, we have the following

\begin{corollary}[Qin \cite{Qin10}, Nakajima \cite{Nakajima09}]
Suppose that $Q$ does not have
oriented cycles. Then all cluster variables of $\ca_Q$ belong to
$\N[x_1^{\pm}, \ldots, x_n^{\pm}]$.
\end{corollary}

This settles a conjecture of Fomin-Zelevinsky \cite{FominZelevinsky02}
in the case of cluster algebras associated with acyclic quivers. The proof
is based on Lusztig's \cite{Lusztig98} and in this sense it does not quite
live up to the hopes that cluster theory ought to explain Lusztig's
results. However, it does show that the conjecture is true for this
important class of cluster algebras.

\section{Caldero-Chapoton's formula}

We describe the bijection
of part a) of theorem~\ref{thm:2}. Let $k$ be an algebraically closed
field and $Q$ a finite quiver without oriented cycles with
vertex set $\{1, \ldots, n\}$. Let $L$ be an object of the
cluster category $\cc_Q$. With $L$, we will associate an element
$X_L$ of the field $\Q(x_1, \ldots, x_n)$. According to
\cite{BuanMarshReinekeReitenTodorov04}, the object $L$
decomposes into a sum of indecomposables
$L_i$, $1\leq i\leq s$, unique up to isomorphism and
permutation. By defining
\[
X_L=\prod_{i=1}^s X_{L_i}
\]
we reduce to the case where $L$ is indecomposable.
Now again by \cite{BuanMarshReinekeReitenTodorov04}, if $L$
is indecomposable, it is either isomorphic to an object $\pi(V)$, or
an object $\Sigma \pi(P_i)$, where $\pi : \cd_Q \to \cc_Q$ is
the canonical projection functor, $\Sigma$ is the suspension functor
of $\cc_Q$, $V$ is a representation of $Q$ (identified with a complex
of representations concentrated in degree $0$) and $P_i$ is the projective
representation associated with a vertex $i$ ($P_i$ is characterized by the existence
of a functorial isomorphism
\[
\Hom(P_i, W) = W_i
\]
for each representation $W$). If $L$ is isomorphic to $\Sigma \pi(P_i)$,
we put $X_L=x_i$. If $L$ is isomorphic to $\pi(V)$, we define
\[
X_L = X_V = \frac{1}{\prod_{i=1}^n x_i^{d_i}} \sum_{0\leq e \leq d} \chi(\mbox{Gr}_e(V))
\prod_{i=1}^n x_i^{\sum_{j\to i} e_j + \sum_{i\to j} (d_j-e_j)} \ko
\]
where $d_i=\dim V_i$, $1 \leq i \leq n$, the sum is taken over all
elements $e\in \N^n$ such that $0 \leq e_i \leq d_i$ for all $i$,
the {\em quiver Grassmannian} $\mbox{Gr}_e(V)$ is the variety
of $n$-tuples of subspaces $U_i \subset V_i$ such that $\dim U_i=e_i$
and the $U_i$ form a subrepresentation of $V$, the Euler characteristic $\chi$
is taken with respect to \'etale cohomology
(or with respect to singular cohomology with coefficients in
a field if $k=\C$) and the sums in the exponent of $x_i$
are taken over all arrows $j \to i$ respectively all
arrows $i\to j$. This formula was invented by P.~Caldero and
F.~Chapoton in \cite{CalderoChapoton06} for the case of
a quiver whose underlying graph is a simply laced Dynkin diagram.
It is still valid for arbitrary quivers without oriented
cycles \cite{CalderoKeller06} and further generalizes to arbitrary
triangulated $2$-Calabi-Yau categories containing a cluster-tilting
object \cite{Palu07}.

\section{Some further developments}

The extension of the results presented here to quivers containing
oriented cycles is the subject of ongoing research.  In a series of
papers \cite{GeissLeclercSchroeer05} \cite{GeissLeclercSchroeer05b} \cite{GeissLeclercSchroeer06}
\cite{GeissLeclercSchroeer06a} \cite{GeissLeclercSchroeer06b}, Geiss-Leclerc-Schr\"oer have obtained
remarkable results for a
class of quivers which are important in the study of (dual
semi-)canonical bases. They use an analogue \cite{GeissLeclercSchroeer05c}
of the Caldero-Chapoton map due ultimately to Lusztig \cite{Lusztig00}.
The class they consider has been further
enlarged by Buan-Iyama-Reiten-Scott \cite{BuanIyamaReitenScott07}.
Thanks to their results, an analogue
of Caldero-Chapoton's formula and a version of theorem~\ref{thm:2}
was proved in \cite{FuKeller07} for an even larger class.

Building on \cite{MarshReinekeZelevinsky03} Derksen-Weyman-Zelevinsky
are developing a represen\-tation-theoretic model for mutation of
general quivers in \cite{DerksenWeymanZelevinsky07}. Their approach
is related to Kontsevich-Soibelman's work \cite{KontsevichSoibelman07},
where $3$-Calabi-Yau categories play an important r\^ole,
as was already the case in \cite{IyamaReiten06}.


\def\cprime{$'$}
\providecommand{\bysame}{\leavevmode\hbox to3em{\hrulefill}\thinspace}
\providecommand{\MR}{\relax\ifhmode\unskip\space\fi MR }
\providecommand{\MRhref}[2]{%
  \href{http://www.ams.org/mathscinet-getitem?mr=#1}{#2}
}
\providecommand{\href}[2]{#2}

\end{document}